\documentclass[11pt]{amsart}
\usepackage{amssymb}
\newcommand{\Z}{\mathbb{Z}}
\newcommand{\N}{\mathbb{N}}

\newcommand{\Ext}{\operatorname{Ext}}

\newcommand{\soc}{\operatorname{soc}}

\newcommand{\Hom}{\operatorname{Hom}}
\newcommand{\PHom}{\operatorname{PHom}}

\newcommand{\End}{\operatorname{End}}

\newcommand{\md}{\operatorname{mod}}
\newcommand{\proj}{\operatorname{proj}}

\newcommand{\sHom}{\operatorname{\underline{Hom}}}
\newcommand{\lra}{\longrightarrow}

\newtheorem{thm}{Theorem}[section]
\newtheorem{lem}[thm]{Lemma}
\newtheorem{prop}[thm]{Proposition}
\newtheorem{cor}[thm]{Corollary}
\newtheorem{conj}[thm]{Conjecture}

\theoremstyle{remark}
\newtheorem{rem}{Remark\!\!}

\theoremstyle{remark}

\theoremstyle{definition}

\newcommand{\colim}{\operatorname{colim}}
\newcommand{\hocolim}{\operatorname{hocolim}}
\newcommand{\Mod}{\operatorname{Mod}}
\newcommand{\stmod}{\operatorname{\sf stmod}}

\numberwithin{equation}{section}

\date{} 
\title{Equivalences of derived categories for symmetric algebras}

\author{Jeremy Rickard}
\address{University of Bristol\\ School of Mathematics\\ 
University Walk\\ Bristol BS8 1TW\\ ENGLAND}
\email{j.rickard@bristol.ac.uk}

\begin{document}

\maketitle

\section{Introduction}
\label{sec:intro}

It is about a decade since Brou\'{e} made his celebrated
conjecture~\cite{bro:perfect} on equivalences of derived categories in
block theory: that the module categories of a block algebra $A$ of a
finite group algebra and its Brauer correspondent $B$ should have
equivalent derived categories if their defect group is abelian. Since
then, character-theoretic evidence for the conjecture has accumulated
rapidly, but until very recently there have been very few examples
where the conjecture has actually been verified. This is because the
precise structure of, say, the indecomposable projective modules for
$A$ is known only in the simplest cases (although the corresponding
structure for $B$ is much easier to determine): this makes it very
difficult to carry out explicit calculations to verify an equivalence
of derived categories.

Recently, however, Okuyama~\cite{oku:derived} introduced a method of
proving that there is an equivalence of derived categories that needs
very little explicit information about $A$. 

In many of the simpler cases where Brou\'{e}'s conjecture is not yet
known to be true, there is known to be a `stable equivalence of Morita
type' between $A$ and $B$: an exact functor between the module
categories that is an equivalence of categories `modulo projective
modules'. This is a consequence of an equivalence of derived
categories, since the stable module category is a canonical quotient
of the derived category. Moreover, recent work of
Rouquier~\cite{rou:gluing},\cite{rou:block} gives a method of
constructing such stable equivalences from equivalences of derived
categories for smaller groups.

Okuyama's method is a strategy for lifting stable equivalences to
equivalences of derived categories. If one can produce an equivalence
of derived categories between $B$ and a third algebra $C$, and if the
objects of the derived category $D^b(\md(B))$ that correspond to the
simple $C$-modules are isomorphic in the stable module category of $B$
(regarded as a quotient category of $D^b(\md(B))$) to the images of
the simple $A$-modules under a stable equivalence of Morita type, then
it follows from a theorem of Linckelmann~\cite[Theorem
2.1]{lin:Morita} that $A$ and $C$ are Morita equivalent, and so $A$
and $B$ have equivalent derived categories. Note that to carry out
this strategy, one needs to know nothing about $A$ except which
objects of the stable module category of $B$ correspond to the simple
$A$-modules.

Okuyama used this method to verify Brou\'{e}'s conjecture for many
examples of blocks with defect group $C_3\times
C_3$~\cite{oku:derived}.

This still leaves the problem of finding a suitable equivalence 
\begin{equation}
  \label{eq:equiv}
  D^b(\md(B))\approx D^b(\md(C))
\end{equation}
of derived categories. Okuyama did this by constructing a suitable
`tilting complex' $T$ for $B$: by the main theorem
of~\cite{ric:Morita}, $T$ is the object corresponding to the free
$C$-module under an equivalence~(\ref{eq:equiv}), where $C$ is the
endomorphism algebra of $T$. The main theorem of this paper,
Theorem~\ref{th:main}, gives an alternative approach. Rather than
characterizing the objects of $D^b(\md(B))$ that correspond to free
modules under equivalences of derived categories, as the definition of
a tilting complex does, we characterize the sets of objects that can
correspond to the simple modules. Our proof requires only that $B$
should be a symmetric algebra, which is of course the case for a block
of a finite group algebra; it is easy to construct counterexamples for
general finite-dimensional algebras, but we do not know any
counterexamples for self-injective algebras.

Since Linckelmann's theorem focuses on the simple modules, this new
characterization of derived equivalence is well suited to applying his
theorem as Okuyama did. In Section~\ref{sec:examples} we give several
fairly simple examples. Simple modules have a less complicated
structure than projective modules, so it is no surprise that in our
examples the objects of $D^b(\md(B))$ that we construct (corresponding
to the simple $A$-modules) are considerably simpler than the tilting
complex (corresponding to a free $A$-module) would be.

In a recent paper~\cite{chu:sl2}, Chuang has used our main theorem to
give a proof of Brou\'{e}'s conjecture for the principal block of
$SL(2,p^2)$ in characteristic $p$ for an arbitrary prime $p$. Using
different methods, this result has subsequently been extended to
$SL(2,p^n)$ for all $n$ by Okuyama~\cite{oku:sl2}.
Holloway~\cite{hol:thesis} has also used our main theorem to verify
Brou\'{e}'s conjecture for several blocks, including the principal
block of the sporadic Hall-Janko group in characteristic five, using
computer calculations.

\section{Conventions and notation}
\label{sec:notation}

Throughout most of this paper, $k$ will be an algebraically closed
field. This is not really essential, but for the sake of clarity we
prefer to avoid the (entirely routine) added complication involved in
dealing with a field that is not algebraically closed.  In
Section~\ref{sec:non-closed} we shall give details of the minor
modifications that need to be made to deal with a general field.

By a `module' for a ring, we shall mean a \emph{left} module unless we
specify otherwise.

If $\Lambda$ is a ring, then $\Mod(\Lambda)$ will be the category of
all (left) $\Lambda$-modules. If $\Lambda$ is a finite-dimensional
$k$-algebra, then $\md(\Lambda)$ will be the category of finitely
generated $\Lambda$-modules and $\proj(\Lambda)$ will be the category
of finitely generated projective $\Lambda$-modules.

If $\Lambda$ is a finite-dimensional self-injective $k$-algebra (i.e.,
the injective and projective $\Lambda$-modules coincide), then
$\stmod(\Lambda)$ will be the stable module category, which is the
quotient of $\md(\Lambda)$ by the ideal of maps that factor through
projective modules. This is a triangulated category with shift functor
the inverse $\Omega^{-1}$ of the Heller translate. The space of maps
from $M$ to $N$ in $\stmod(\Lambda)$ will be denoted by
$\sHom_{\Lambda}(M,N)$.

Our `complexes' will all be \emph{cochain} complexes, so the
differentials will have degree +1.

If
$$X:= \dots\longrightarrow
X^{i-1}\stackrel{d_X^{i-1}}{\longrightarrow}X^i
\stackrel{d_X^i}{\longrightarrow}X^{i+1}\longrightarrow\dots$$
is a
cochain complex, then $X[m]$ will be $X$ `shifted $m$ places to the
left': i.e., $X[m]^i=X^{i+m}$ and $d_{X[m]}^i=(-1)^md_X^{i+m}$.

If $\mathcal{A}$ is an additive category, then $K(\mathcal{A})$ will
be the chain homotopy category of cochain complexes over $\mathcal{A}$,
$K^-(\mathcal{A})$ will be the full subcategory consisting of
complexes $X$ that are `bounded above' (i.e., $X^i=0$ for $i>>0$) and
$K^b(\mathcal{A})$ will be the full subcategory of `bounded' complexes
(i.e., complexes $X$ with $X^i=0$ for all but finitely many $i$).

If $\mathcal{A}$ is an abelian category, then $D(\mathcal{A})$ will be
the derived category of cochain complexes over $\mathcal{A}$, and
$D^-(\mathcal{A})$ and $D^b(\mathcal{A})$ will be the full
subcategories of complexes that are respectively bounded above and
bounded.

We shall regard an abelian category $\mathcal{A}$ as a full
subcategory of its derived category $D(\mathcal{A})$ in the usual way,
identifying an object $X$ of $\mathcal{A}$ with the complex whose only
non-zero term is $X$ in degree zero.

\section{Preliminaries on symmetric algebras}
\label{sec:symmetric}

A finite-dimensional $k$-algebra $\Lambda$ is said to be
\textbf{symmetric} if there is a \textbf{symmetrizing form} on
$\Lambda$: i.e., a linear map $\lambda:\Lambda\rightarrow k$ such that
$$\forall x,y\in\Lambda, \lambda(xy)=\lambda(yx),$$
and such that the
kernel of $\lambda$ contains no non-zero left or right ideal of
$\Lambda$.

The principal example is the group algebra $kG$ of a finite group,
when $\lambda$ can be taken to be
$$\lambda(g)=\Bigl\{
\begin{array}{ll}1 & \textrm{if $g=1$} \\
0 & \textrm{if $g\in G-\{1\}$,}
\end{array}$$
extended linearly to the whole of $kG$. Any block algebra of $kG$ is
also a symmetric algebra, using the restriction of the same map
$\lambda$ to the block.

The following theorem, giving characterizations of symmetric algebras
in terms of the module category, is well-known, but for the readers'
convenience we include a proof.

\begin{thm}
\label{th:symmetric}
Let $\Lambda$ be a finite-dimensional $k$-algebra. The following
conditions are equivalent.
  \begin{itemize}
  \item[$(a)$] $\Lambda$ is symmetric.
  \item[$(b)$] $\Lambda$ and its dual
    $\Lambda^{\vee}=\Hom_k(\Lambda,k)$ are isomorphic as
    $\Lambda$-bimodules.
  \item[$(c)$] $\Hom_k(?,k)$ and $\Hom_{\Lambda}(?,\Lambda)$ are
    isomorphic as functors from the category of left $\Lambda$-modules
    to the category of right $\Lambda$-modules.
  \item[$(c')$] $\Hom_k(?,k)$ and $\Hom_{\Lambda}(?,\Lambda)$ are
    isomorphic as functors from the category of right
    $\Lambda$-modules to the category of left $\Lambda$-modules.
  \item[$(d)$] For finitely generated projective left
    $\Lambda$-modules $P$ and finitely generated left
    $\Lambda$-modules $M$, there is an isomorphism of $k$-vector
    spaces
    $$\Hom_{\Lambda}(P,M)\cong\Hom_{\Lambda}(M,P)^{\vee},$$
    functorial
    in both $P$ and $M$.
  \item[$(d')$] For finitely generated projective right
    $\Lambda$-modules $P$ and finitely generated right
    $\Lambda$-modules $M$, there is an isomorphism of $k$-vector
    spaces
    $$\Hom_{\Lambda}(P,M)\cong\Hom_{\Lambda}(M,P)^{\vee},$$
    functorial
    in both $P$ and $M$.
  \end{itemize}
\end{thm}

\begin{proof}
  Since (a) and (b) are left-right symmetric, it is sufficient to
  prove that (a) and (b) are equivalent and that
  $\textrm{(b)}\Rightarrow\textrm{(c)}\Rightarrow
  \textrm{(d)}\Rightarrow\textrm{(b)}$.
  
  First we shall show that (a) implies (b). Let
  $\lambda:\Lambda\rightarrow k$ be a symmetrizing form on $\Lambda$.
  Consider the map
  $$\theta:\Lambda\rightarrow\Hom_k(\Lambda,k)$$
  defined by
  $\theta(x)=x.\lambda$ for $x\in\Lambda$: i.e.,
  $\theta(x):\Lambda\rightarrow k$ is the map
  $$y\mapsto\lambda(yx)=\lambda(xy).$$
  Since $\theta(xz)$ is the map
  $$y\mapsto\lambda(yxz)=\lambda(zyx) =(x.\lambda.z)(y),$$
  $\theta$ is
  a $\Lambda$-bimodule homomorphism.
  
  Now let $x\in\Lambda$, and suppose $\theta(x)=0$. Then
  $\lambda(yx)=0$ for all $y\in\Lambda$, and so the left ideal
  $\Lambda x$ is contained in the kernel of $\lambda$. Since $\lambda$
  is a symmetrizing form, this implies $x=0$. Hence $\theta$ is
  injective, and is therefore an isomorphism, since $\Lambda$ and
  $\Lambda^{\vee}$ are vector spaces of the same dimension.
  
  To show (b) implies (a), suppose
  $\theta:\Lambda\rightarrow\Lambda^{\vee}$ is a $\Lambda$-bimodule
  isomorphism. Let $\lambda=\theta(1)$. Then for $x,y\in\Lambda$,
  $$\lambda(xy)=\theta(1)(xy)= (y\theta(1))(x)=\theta(y)(x)=
  (\theta(1)y)(x)=\theta(1)(yx)= \lambda(yx).$$
  Also, suppose the left
  ideal $\Lambda x$ is contained in the kernel of $\lambda$. Then for
  every $y\in\Lambda$,
  
  $$0=\lambda(yx)=(x\theta(1))(y)= \theta(x)(y),$$
  and so
  $\theta(x)=0$, and so $x=0$, since $\theta$ is an isomorphism.
  Similarly, no non-zero right ideal of $\Lambda$ is contained in the
  kernel of $\lambda$, and so $\lambda$ is a symmetrizing form on
  $\Lambda$.
  
  Next we show that (b) implies (c). Suppose that $\Lambda$ and
  $\Lambda^{\vee}$ are isomorphic as $\Lambda$-bimodules. Let $M$ be a
  left $\Lambda$-module. There is a chain of natural isomorphisms of
  right $\Lambda$-modules
  $$\Hom_{\Lambda}(M,\Lambda)\cong\Hom_{\Lambda}(M,\Hom_k(\Lambda,k))
  \cong\Hom_k(\Lambda\otimes_{\Lambda}M,k)\cong\Hom_k(M,k).$$
  
  Next assume that (c) is true, and let us deduce (d). Let $M$ and $P$
  be finitely generated left $\Lambda$-modules. Then
  $\Hom_{\Lambda}(M,P)\cong\Hom_k(P^{\vee}\otimes_{\Lambda}M,k)$, and
  since these are all finite-dimensional vector spaces,
  $\Hom_{\Lambda}(M,P)^{\vee}\cong P^{\vee}\otimes_{\Lambda}M$, which,
  since (c) is true, is in turn isomorphic to
  $\Hom_{\Lambda}(P,\Lambda)\otimes_{\Lambda}M$. There is a natural
  map
  \begin{equation}
    \label{eq:symm}
    \Hom_{\Lambda}(P,\Lambda)\otimes_{\Lambda}M\longrightarrow
    \Hom_{\Lambda}(P,M),
  \end{equation}
  sending $\alpha\otimes m$ (where
  $\alpha\in\Hom_{\Lambda}(P,\Lambda)$ and $m\in M$) to the map
  $P\rightarrow M$ sending $p\in P$ to $\alpha(p)m$. For $P=\Lambda$,
  it is easy to check that (\ref{eq:symm}) is an isomorphism between
  vector spaces naturally isomorphic to $M$. By naturality,
  (\ref{eq:symm}) is also an isomorphism for $P$ any direct summand of
  a finite direct sum of copies of $\Lambda$; i.e., for any finitely
  generated projective module $P$.
  
  Finally, the fact that (d) implies (b) follows by taking
  $M=P={}_{\Lambda}\Lambda$. The natural isomorphism
  $$\Lambda\cong \Hom_{\Lambda}(\Lambda,\Lambda)$$
  is an isomorphism
  of $\Lambda$-bimodules, where the bimodule structure on the right
  hand side is induced by the \emph{right} action of $\Lambda$ by
  multiplication on the two arguments. But by (d) the right hand side
  is naturally isomorphic to its dual, and by naturality this
  isomorphism is an isomorphism of bimodules. Hence (b).
\end{proof}

We shall most often be using condition (d).

Here is a corollary about maps in the derived category.

\begin{cor}
  \label{cor:symm}
  Let $\Lambda$ be a finite-dimensional symmetric $k$-algebra. Let
  $P^*$ be a bounded complex of finitely generated projective
  left $\Lambda$-modules, and let $M^*$ be a complex of finitely
  generated left $\Lambda$-modules. Then
  $\Hom_{D(\Mod(\Lambda))}(P^*,M^*)$ and
  $\Hom_{D(\Mod(\Lambda))}(M^*,P^*)$ are naturally dual.
\end{cor}

\begin{proof}
  Since $P^*$ is a bounded complex of projective and (since $\Lambda$
  is symmetric) injective modules, calculating homomorphisms from or
  to $P^*$ in the derived category is equivalent to doing so in the
  chain homotopy category $K(\Mod(\Lambda))$.
  
  Recall that if $X^*$ and $Y^*$ are complexes of
  $\Lambda$-modules, then $\Hom_{K(\Mod(\Lambda))}(X^*,Y^*)$ may be
  calculated as the degree zero homology of the `completed' total
  complex of the double complex $\Hom_{\Lambda}(X^*,Y^*)$ (i.e., the
  variation of the total complex where the terms are formed by taking
  direct \emph{products} rather than the direct \emph{sums} of
  diagonals in the double complex). In fact, the fact that $P^*$ is
  bounded ensures that for $X^*=P^*$ and $Y^*=M^*$ or vice versa,
  there is only a finite number of non-zero terms on each diagonal,
  and so the completed total complex is the same as the usual total
  complex.
  
  Condition (d) of Theorem~\ref{th:symmetric} implies that the double
  complexes $\Hom_{\Lambda}(P^*,M^*)$ and $\Hom_{\Lambda}(M^*,P^*)$
  are naturally dual. By the remark at the end of the last paragraph,
  their total complexes are complexes of finite dimensional vector
  spaces and are also naturally dual. Taking degree zero homology, the
  corollary follows.
\end{proof}

\section{Preliminaries on homotopy colimits}
\label{sec:hocolim}

In this section, $\Lambda$ will be an arbitrary ring.

Let
$$X_0\stackrel{\alpha_0}{\lra} X_1 \stackrel{\alpha_1}{\lra} X_2
\stackrel{\alpha_2}{\lra}\dots$$
be a sequence of maps in a
triangulated category $\mathcal{T}$ with countable coproducts.
Recall~\cite{BN:hocolim} that the \textbf{homotopy colimit}
$\hocolim(X_i)$ of this sequence is defined by forming the
distinguished triangle
\begin{equation}
  \label{eq:hocolim}
  \bigoplus_{i=0}^{\infty}X_i\lra\bigoplus_{i=0}^{\infty}X_i\lra
  \hocolim(X_i)\lra\bigoplus_{i=0}^{\infty}X_i[1],
\end{equation}
where the restriction of the first map to $X_i$ is
$\mathrm{id}_{X_i}-\alpha_i$. This defines $\hocolim(X_i)$ up to
isomorphism, but not usually up to \emph{unique} isomorphism.

If $\mathcal{T}$ is the derived category $D(\Mod(\Lambda))$, and if we
can choose chain maps $\beta_i$ representing the maps $\alpha_i$ of
the derived category, then $\hocolim(X_i)$ is isomorphic to the usual
colimit in the category of chain complexes of the sequence
$$X_0\stackrel{\beta_0}{\lra} X_1 \stackrel{\beta_1}{\lra} X_2
\stackrel{\beta_2}{\lra}\dots.$$
This is an easy consequence of the
fact that the coproduct, in the derived category, of a family of
complexes is the same as the coproduct in the category of chain
complexes.

An object $C$ of the triangulated category $\mathcal{T}$ is called
\textbf{compact} if the functor $\Hom(C,?)$ commutes with arbitrary
coproducts: more precisely, if the natural map
$$\bigoplus_{i\in I}\Hom(C,X_i)\lra \Hom\bigl(C,\bigoplus_{i\in
  I}X_i)$$
is an isomorphism whenever $\{X_i:i\in I\}$ is a set of
objects of $\mathcal{T}$ whose coproduct exists in $\mathcal{T}$.

In the derived category $D(\Mod(\Lambda))$ of a module category, the
compact objects are precisely those that are isomorphic to bounded
complexes of finitely generated projective modules.

An easy consequence of the definition is that if $C$ is a compact
object, then there is a natural isomorphism
$$\Hom\bigl(C,\hocolim(X_i)\bigr)\lra \colim\bigl(\Hom(C,X_i)\bigr).$$

We shall need a generalization of this. Let us say that a family
$\mathcal{X}$ of objects of the derived category $D(\Mod(\Lambda))$ is
\textbf{uniformly bounded below} if there is some $n\in\Z$ such that
the degree $j$ cohomology $H^j(X)$ of $X$ is zero for all
$X\in\mathcal{X}$ and all $j<n$.

\begin{prop}
  \label{prop:uniform}
  Let $C$ be an object of $D(\Mod(\Lambda))$ isomorphic to a complex
  of finitely generated projective $\Lambda$-modules that is bounded
  above. For example, if $\Lambda$ is a finite-dimensional
  $k$-algebra, let $C$ be any object of $D^-(\md(\Lambda))$.
  
  (a) Let $\mathcal{X}$ be a family of objects of $D(\Mod(\Lambda))$
  that is uniformly bounded below. Then the natural map
  
  $$\bigoplus_{X\in\mathcal{X}}\Hom(C,X)\longrightarrow
  \Hom\bigl(C,\bigoplus_{X\in\mathcal{X}}X\bigr)$$
  is an isomorphism.
  
  (b) Let $X_0\rightarrow X_1\rightarrow X_2\rightarrow\dots$ be a
  sequence of maps in the derived category $D(\Mod(\Lambda))$, and
  suppose that $\{X_i:i\in\N\}$ is uniformly bounded below. Then
  $$\Hom\bigl(C,\hocolim(X_i)\bigr)\cong
  \colim\bigl(\Hom(C,X_i)\bigr).$$
\end{prop}

\begin{proof}
  First, if $\Lambda$ is a finite-dimensional $k$-algebra, then any
  object of $D^-(\md(\Lambda))$ is isomorphic to its minimal
  projective resolution, which is a complex of finitely generated
  projective $\Lambda$-modules that is bounded above.
  
  (a) Without loss of generality, we shall assume that every
  $X\in\mathcal{X}$ is a complex with no non-zero terms in negative
  degrees, and that $C$ is the bounded above complex
  $$\dots\lra C^{-2}\lra C^{-1}\lra C^0\lra C^1\lra\dots,$$
  of finitely
  generated projectives. Let $\tilde{C}$ be the truncated complex
  $$\dots\lra 0\lra C^{-1}\lra C^0\lra C^1\lra\dots,$$
  which is a
  bounded complex of finitely generated projectives.  The inclusion
  $\tilde{C}\rightarrow C$ of complexes induces an isomorphism
  $$\Hom_{D(\Mod(\Lambda))}(C,Y)\longrightarrow
  \Hom_{D(\Mod(\Lambda))}(\tilde{C},Y)$$
  for $Y\in\mathcal{X}$ and for
  $Y=\bigoplus_{X\in\mathcal{X}}X$. The proposition follows because
  $\tilde{C}$ is compact.
  
  (b) This follows because, by (a), the natural map
  $$\bigoplus_{i\in\N}\Hom(C,X_i[m])\longrightarrow
  \Hom\bigl(C,\bigoplus_{i\in\N}X_i[m]\bigr)$$
  is an isomorphism for
  every $m$, and so the long exact sequence obtained by applying the
  functor $\Hom(C,?)$ to the triangle~(\ref{eq:hocolim}) breaks up
  into a sequence of short exact sequences, including one isomorphic
  to
  $$0\longrightarrow\bigoplus_{i\in\N}\Hom(C,X_i) \longrightarrow
  \bigoplus_{i\in\N}\Hom(C,X_i)\longrightarrow
  \Hom(C,\hocolim(X_i))\longrightarrow0,$$
  which expresses
  $\Hom(C,\hocolim(X_i))$ as the colimit of $\Hom(C,X_i)$.
\end{proof}

\section{The main theorem}
\label{sec:main}

Let $\Lambda$ and $\Gamma$ be two rings. A necessary and sufficient
condition~\cite{ric:Morita} for the derived categories
$D(\Mod(\Lambda))$ and $D(\Mod(\Gamma))$ to be equivalent as
triangulated categories is that there should be a \textbf{tilting
  complex} in $D(\Mod(\Lambda))$ whose endomorphism ring is isomorphic
to $\Gamma$: i.e., a bounded complex $T$ of finitely generated
projective $\Lambda$-modules such that
\begin{itemize}
\item[(i)] $\Hom(T,T[m])=0$ for $m\neq0$, and
\item[(ii)] the direct summands of $T$ generate $K^b(\proj(\Lambda))$
  as a triangulated category,
\end{itemize}
with $\End(T)\cong\Gamma$.

If $T$ is such a tilting complex, then there is an equivalence
$D(\Mod(\Lambda))\approx D(\Mod(\Gamma))$ that sends $T$ to the free
$\Gamma$-module of rank one, so the indecomposable summands of $T$
correspond to the indecomposable projective $\Gamma$-modules.

We shall use later the fact~\cite{ric:Morita} that, for a
finite-dimensional algebra $\Lambda$, condition (ii) can be replaced
by the condition that, for any non-zero object $C$ of
$D^-(\md(\Lambda))$, $\Hom(T,C[m])\neq0$ for some $m\in\Z$.

In this section we shall consider instead the objects $X_0,\dots,X_r$
of $D(\Mod(\Lambda))$ that correspond to the simple $\Gamma$-modules
in the case that $\Lambda$ and $\Gamma$ are finite-dimensional
symmetric $k$-algebras.  Since an equivalence $D(\Mod(\Lambda))\approx
D(\Mod(\Gamma))$ restricts to an equivalence between the full
subcategories of objects isomorphic to bounded complexes of finitely
generated modules, $X_0,\dots,X_r$ must (up to isomorphism) be objects
of $D^b(\md(\Lambda))$. Also, for $0\leq i,j\leq r$, they must satisfy
\begin{itemize}
\item[(a)] $\Hom(X_i,X_j[m])=0$ for $m<0$,
\item[(b)] $\Hom(X_i,X_j)=\Bigl\{
    \begin{array}{ll}
      0 & \textrm{if } i\neq j \\
      k & \textrm{if } i=j
    \end{array}$, and
  \item[(c)] $X_0,\dots,X_r$ generate $D^b(\md(\Lambda))$ as a
    triangulated category,
\end{itemize}
since the simple $\Gamma$-modules satisfy the corresponding
properties.

In this section we shall prove a partial converse to this.

\begin{thm}
  \label{th:main}
  Let $\Lambda$ be a finite-dimensional symmetric $k$-algebra.  Let
  $X_0,\dots,X_r$ be objects of $D^b(\md(\Lambda))$ satisfying
  conditions (a)--(c) above. Then there is another finite-dimensional
  symmetric $k$-algebra $\Gamma$ and an equivalence of triangulated
  categories
  $$D(\Mod(\Lambda))\approx D(\Mod(\Gamma))$$
  sending $X_0,\dots,X_r$
  to the simple $\Gamma$-modules.
\end{thm}

We shall give the proof as a sequence of lemmas. 

What we shall do is construct a tilting complex
$T=T_0\oplus\dots\oplus T_r$ for $\Lambda$ such that, for $0\leq
i,j\leq r$ and $m\in\Z$,
\begin{equation}
  \label{eq:TtoX}
  \Hom(T_i,X_j[m])=\Bigl\{
  \begin{array}{ll}
    k & \textrm{if } i=j \textrm{ and } m=0 \\
    0 & \textrm{otherwise.}
  \end{array}
\end{equation}

\begin{lem}
  \label{lem:reduction}
  Suppose there is a tilting complex $T=T_0\oplus\dots\oplus T_r$ for
  $\Lambda$ satisfying the property (\ref{eq:TtoX}). Then
    Theorem~\ref{th:main} is true.
\end{lem}

\begin{proof}
  For $\Gamma=\End(T)$, there is an equivalence of
  triangulated categories
  $$F:D(\Mod(\Lambda))\rightarrow D(\Mod(\Gamma))$$
  sending
  $T_0,\dots,T_r$ to the indecomposable projective $\Gamma$-modules
  $P_0,\dots,P_r$. Since the simple $\Gamma$-modules $S_0,\dots,S_r$,
  numbered so that $P_i$ is the projective cover of $S_i$, are
  characterized up to isomorphism in $D(\Mod(\Gamma))$ by the fact
  that, for $0\leq i,j\leq r$ and $m\in\Z$,
  $$\Hom(P_i,S_j[m])=\Bigl\{
  \begin{array}{ll}
    k & \textrm{if } i=j \textrm{ and } m=0 \\
    0 & \textrm{otherwise}
  \end{array},$$
  the equivalence $F$ sends $X_0,\dots,X_r$ to the simple modules
  $S_0,\dots,S_r$. 
  
  Since a ring whose derived category is equivalent to a
  finite-dimensional symmetric $k$-algebra is itself a
  finite-dimensional symmetric $k$-algebra~\cite{ric:equivalences},
  Theorem~\ref{th:main} follows.
\end{proof}
  
Now let us construct the summands $T_i$ of the tilting complex $T$. 

Set $X_i^{(0)}:=X_i$. By induction on $n$, we shall construct a
sequence 
$$X_i^{(0)}\rightarrow X_i^{(1)}\rightarrow\dots\rightarrow
X_i^{(n-1)}\rightarrow X_i^{(n)}\rightarrow\dots$$
of objects and maps in $D(\Mod(\Lambda))$.

Suppose we have constructed $X_i^{(n-1)}$.  For each $0\leq j\leq r$
and $t<0$, choose a basis $B_i^{(n-1)}(j,t)$ of
$\Hom(X_j[t],X_i^{(n-1)})$, let $Z_i^{(n-1)}(j,t)$
be a direct sum of copies of $X_j[t]$ indexed by $B_i^{(n-1)}(j,t)$,
and let
$$\alpha_i^{(n-1)}(j,t): Z_i^{(n-1)}(j,t)\rightarrow X_i^{(n-1)}$$
be the map
whose restriction to the summand of $Z_i^{(n-1)}(j,t)$ corresponding
to an element $\beta\in B_i^{(n-1)}(j,t)$ is just $\beta$. Now let 
$$Z_i^{(n-1)}= \mathop{\bigoplus_{0\leq j\leq
    r}}_{t<0}Z_i^{(n-1)}(j,t),$$
and let
$$\alpha_i^{(n-1)}: Z_i^{(n-1)}\rightarrow X_i^{(n-1)}$$
be the map whose restriction to the summand $Z_i^{(n-1)}(j,t)$ is
$\alpha_i^{(n-1)}(j,t)$. 

All that will matter about this map is that
\begin{itemize}
\item[(a)] $Z_i^{(n-1)}$ is a (possibly infinite) direct sum of
  objects, each of which is isomorphic to $X_j[t]$ for some $0\leq
  j\leq r$ and some $t<0$,
\item[(b)] the map $\Hom(X_j[t], Z_i^{(n-1)})
  \rightarrow\Hom(X_j[t], X_i^{(n-1)})$ induced by
  $\alpha_i^{(n-1)}$ is surjective for every $0\leq j\leq r$ and
  $t<0$, and
\item[(c)] the map $\Hom(X_j[-1], Z_i^{(n-1)})
  \rightarrow\Hom(X_j[-1], X_i^{(n-1)})$ induced by
  $\alpha_i^{(n-1)}$ is an isomorphism for every $0\leq j\leq r$.
\end{itemize}

The map we have constructed clearly satisfies properties (a) and (b). Since
$$\Hom(X_j[-1], X_{j'}[t])=
\Hom(X_j,X_{j'}[t+1])=0$$
for all $0\leq j,j'\leq
r$ and $t<0$ except for $j=j'$ and $t=-1$, property (c) is also
satisfied.

Now define $X_i^{(n)}$, together with a map from $X_i^{(n-1)}$, by
forming the distinguished triangle
\begin{equation}
  \label{eq:triangle}
  Z_i^{(n-1)}\stackrel{\alpha_i^{(n-1)}}{\longrightarrow}X_i^{(n-1)}
  \longrightarrow X_i^{(n)}\longrightarrow Z_i^{(n-1)}[1].
\end{equation}

Finally, define $T_i:=\hocolim(X_i^{(n)})$.

\begin{rem}
  In general it will not be possible to choose $Z_i^{(n-1)}$ in 
  the bounded derived category
  $D^b(\md(\Lambda))$ that satisfy properties (a)--(c), even if
  $X_i^{(n-1)}$ is in $D^b(\md(\Lambda))$. However, if $\Lambda$ is
  the group algebra of a finite group, then one can prove, using the
  fact that the cohomology algebra is finitely generated, that we can
  do so, and so, by induction on $n$, we may force every $X_i^{(n)}$
  to be in $D^b(\md(\Lambda))$. However, it is still not obvious
  \emph{a priori} that the homotopy colimit $T_i$ will be in
  $D^b(\md(\Lambda))$.
\end{rem}

\begin{lem}
  \label{lem:YtoT}
  If $Y$ is an object of $D^-(\md(\Lambda))$, then 
  $$\Hom(Y,T_i)\cong
  \colim\bigl(\Hom(Y,X_i^{(n)})\bigr).$$
\end{lem}

\begin{proof}
  A projective resolution of $Y$ is a complex of finitely generated
  projective $\Lambda$-modules that is bounded above, so the Lemma
  will follow from Proposition~\ref{prop:uniform}(b) if we can show that
  $\{X_i^{(n)}:n\geq0\}$ is uniformly bounded below.

  Since $\{X_i:0\leq i\leq r\}$ is a finite set of objects of
  $D^b(\md(\Lambda))$, we shall assume, without loss of generality,
  that the cohomology $H^m(X_i)$ vanishes for all $m<0$. But then
  $H^m(Z_i^{(n)})=0$ for all $i,n$ and all $m<1$, and so, by induction
  on $n$ and the long exact sequence of homology for the distinguished
  triangle (\ref{eq:triangle}), we have $H^m(X_i^{(n)})=0$ for all
  $m<0$ and all $n$. In other words, $\{X_i^{(n)}:n\geq0\}$ is
  uniformly bounded below, as required.
\end{proof}

\begin{lem}
  \label{lem:XtoT}
  For $0\leq i,j\leq r$, and $m\in\Z$,
  $$\Hom(X_j,T_i[m])=\Bigl\{
  \begin{array}{ll}
    k & \textrm{if } i=j \textrm{ and } m=0 \\
    0 & \textrm{otherwise.}
  \end{array}$$
\end{lem}

\begin{proof}
  We shall prove this lemma by considering the natural map
  \begin{equation}
    \label{eq:XtoX}
    \Hom(X_j,X_i^{(n-1)}[m])\longrightarrow
    \Hom(X_j,X_i^{(n)}[m])
  \end{equation}
  by using the long exact sequence obtained by applying the functor
  $\Hom(X_j,?)$ to the distinguished triangle (\ref{eq:triangle}).
  
  Since, by Proposition~\ref{prop:uniform}(a),
  $\Hom(X_j,Z_i^{(n-1)}[m])=0$ for $m\leq0$, the map (\ref{eq:XtoX})
  is an isomorphism for $m<0$ and is injective for $m=0$.
  
  By property (c) of $\alpha_i^{(n-1)}:Z_i^{(n-1)}\rightarrow
  X_i^{(n-1)}$, the map (\ref{eq:XtoX}) is also surjective for $m=0$.

  Since, for $m\leq0$, 
  $$\Hom(X_j,X_i^{(0)}[m])=\Bigl\{
    \begin{array}{ll}
      k & \textrm{if } i=j \textrm{ and } m=0 \\
      0 & \textrm{otherwise}
    \end{array},$$
    Lemma~\ref{lem:YtoT} implies the cases of the lemma involving
    $m\leq0$.
  
  By property (b) of $\alpha_i^{(n-1)}$, the map (\ref{eq:XtoX}) is
  zero for $m>0$. So Lemma~\ref{lem:YtoT} implies that
  $\Hom(X_j,T_i[m])=0$ for $m>0$.
\end{proof}

\begin{lem}
  \label{lem:compact}
  For each $0\leq i\leq r$, $T_i$ is a compact object of
  $D(\Mod(\Lambda))$: i.e., it is isomorphic to a bounded complex of
  finitely generated projectives.
\end{lem}

\begin{proof}
  For an object $Y$ of $D(\Mod(\Lambda))$, consider the vector space
  $$\bigoplus_{m\in\Z}\Hom(Y,T_i[m]).$$
  By Lemma~\ref{lem:XtoT} this
  is finite-dimensional for $Y=X_j$, for any $0\leq j\leq r$. The
  class of objects $Y$ for which it is finite-dimensional form a full
  triangulated subcategory of $D(\Mod(\Lambda))$, and so, since
  $\{X_j:0\leq j\leq r\}$ generates $D^b(\md(\Lambda))$ as a
  triangulated category, it is finite-dimensional for any $Y$ in
  $D^b(\md(\Lambda))$. 

  In particular, 
  $$\bigoplus_{m\in\Z}H^m(T_i)\cong
  \bigoplus_{m\in\Z}\Hom(\Lambda,T_i[m])$$
  is
  finite-dimensional, so $T_i$ is isomorphic to an object of
  $D^b(\md(\Lambda))$. Then, for any simple module $S$,
  $\bigoplus_{m\in\Z}\Hom(S,T_i[m])$ is
  finite-dimensional, so a minimal injective resolution of $T_i$
  contains only a finite number of copies of the injective hull of
  $S$. Thus this injective resolution is a bounded complex of finitely
  generated injectives, which are also projective, since $\Lambda$ is
  symmetric.
\end{proof}

\begin{lem}
  \label{lem:TtoX}
  For $0\leq i,j\leq r$ and $m\in\Z$, 
  $$\Hom(T_i,X_j[m])=\Bigl\{
  \begin{array}{ll}
    k & \textrm{if } i=j \textrm{ and } m=0 \\
    0 & \textrm{otherwise.}
  \end{array}$$
\end{lem}

\begin{proof}
  Since, by Lemma~\ref{lem:compact}, we now know that $T_i$ is
  isomorphic to a bounded complex of projectives, this follows by
  combining Lemma~\ref{lem:XtoT} and Corollary~\ref{cor:symm}.
\end{proof}

\begin{lem}
  \label{lem:TtoT}
  For $0\leq i,j\leq r$ and $m\neq0$,
  $$\Hom(T_i,T_j[m])=0.$$
\end{lem}

\begin{proof}
  By Lemma~\ref{lem:TtoX}, $\Hom(T_i,X_l[m])=0$ for
  all $m<0$ and $0\leq l\leq r$. Hence, for any $n$,
  $\Hom(T_i,Z_j^{(n)}[m])=0$ for $m\leq0$, since
  $T_i$ is compact (by Lemma~\ref{lem:compact}) and $Z_j^{(n)}$ is a
  direct sum of copies of \emph{negative} shifts of copies of various
  $X_l$s. Applying the functor $\Hom(T_i,?)$ to the
  triangle~\ref{eq:triangle}, we get an exact sequence
  $$\Hom(T_i,X_j^{(n-1)}[m])\rightarrow
  \Hom(T_i,X_j^{(n)}[m])\rightarrow
  \Hom(T_i,Z_j^{(n-1)}[m+1])=0$$
  for every $m<0$.
  So by induction on $n$, $\Hom(T_i,X_j^{(n)}[m])=0$ if
  $m<0$.
  
  By Proposition~\ref{prop:uniform}(b),
  $\Hom(T_i,T_j[m])=0$ if $m<0$. For $m>0$ it
  follows that $\Hom(T_j,T_i[m])=0$ by
  Corollary~\ref{cor:symm}.
\end{proof}

\begin{lem}
  \label{lem:Tperp=0}
  Let $C$ be an object of $D^-(\md(\Lambda))$. If $C\not\cong0$, then
  for some $0\leq i\leq r$ and some $m\in\Z$,
  $\Hom(C,T_i[m])\neq0$.
\end{lem}

\begin{proof}
  Since $C$ is bounded above and $X_i$ is bounded,
  $\Hom(C,X_i[m])=0$ for $m<<0$.
  
  However, if $\Hom(C,X_i[m])=0$ for \emph{all} $i$
  and $m$, then $\Hom(C,X)=0$ for every object $X$
  of $D^b(\md(\Lambda))$, since $\{X_i:0\leq i\leq r\}$ generates
  $D^b(\md(\Lambda))$ as a triangulated category. In particular,
  $\Hom(C,\Lambda[m])=0$ for all $m$, and so
  $C\cong0$.
  
  Thus, if $C\not\cong0$, we can choose $m$ and $i$ so that
  $\Hom(C,X_i[m])\neq0$ and $m$ is minimal: i.e.,
  $\Hom(C,X_j[m'])=0$ for all $0\leq j\leq r$ and
  all $m'<m$.

  Let us apply the functor $\Hom(C,?)$ to the
  triangle~(\ref{eq:triangle}). 
  Since, by Proposition~\ref{prop:uniform}(a),
  $\Hom(C,Z_i^{(n-1)}[m])=0$, we get an exact sequence
  $$0=\Hom(C,Z_i^{(n-1)}[m])\rightarrow
  \Hom(C,X_i^{(n-1)}[m])\rightarrow
  \Hom(C,X_i^{(n)}[m]).$$
  In other words,
  $$\Hom(C,X_i^{(n-1)}[m])\rightarrow
  \Hom(C,X_i^{(n)}[m])$$
  is injective for every
  $n\geq1$. Hence
  $\colim\bigl(\Hom(C,X_i^{(n)}[m])\bigr)\neq0$,
  since $\Hom(C,X_i^{(0)}[m])\neq0$. So by
  Lemma~\ref{lem:YtoT}, $\Hom(C,T_i[m])\neq0$.
\end{proof}

We now have all the ingredients to complete the proof of our main
theorem. 

\begin{proof}[Proof of Theorem~\ref{th:main}]
  By Lemmas~\ref{lem:reduction} and~\ref{lem:TtoX}, it is sufficient
  to show that $T=\bigoplus_{0\leq i\leq r}T_i$ is a tilting complex
  for $\Lambda$.
  
  Lemma~\ref{lem:compact} shows that $T$ is isomorphic in
  $D(\Mod(\Lambda))$ to a bounded complex of finitely generated
  projective modules.
  
  Lemma~\ref{lem:TtoT} shows that $\Hom(T,T[m])=0$
  for $m\neq0$.
  
  We just need to show that if $C$ is an object of $D^-(\md(\Lambda))$
  such that $\Hom(T,C[m])=0$ for all $m\in\Z$, then
  $C\cong0$. Let $C$ be such an object. By Corollary~\ref{cor:symm},
  $$\Hom(C,T[-m])\cong
  \Hom(T,C[m])^{\vee}\cong0$$
  for all $m\in\Z$, and
  so $C\cong0$ by Lemma~\ref{lem:Tperp=0}.
\end{proof}

\section{Lifting stable equivalences}
\label{sec:stable}

Suppose $\Gamma'$ is another algebra satisfying the conclusion of
Theorem~\ref{th:main}. Then we have an equivalence of derived
categories
$$D^b(\md(\Gamma))\approx D^b(\md(\Gamma'))$$
that takes the simple
$\Gamma$-modules to the simple $\Gamma'$-modules. Since the finitely
generated projective modules $P$ are characterized up to isomorphism
in $D^b(\md(\Gamma))$ by the fact that
$\Hom(P,S[m])=0$ for every simple module $S$ and
every integer $m\neq0$, this equivalence of derived categories
restricts to give an equivalence
$$\proj(\Gamma)\approx\proj(\Gamma').$$
Hence the algebra $\Gamma$ is
determined up to Morita equivalence by the objects $X_0,\dots,X_r$.

In general, however, the algebra $\Gamma$ might be hard to identify.
In this section we shall build on an idea of Okuyama to show how this
algebra can be identified in certain cases.

First we shall briefly describe Okuyama's method~\cite{oku:derived} of
lifting stable equivalences to equivalences of derived categories.

\subsection{Okuyama's method}
\label{sec:okuyama}

Let $A$ and $B$ be finite-dimensional self-injective $k$-algebras.

Recall that the stable module category $\stmod(A)$ of $A$ is
equivalent to the quotient of triangulated categories
$$D^b(\md(A))/K^b(\proj(A)).$$
Thus any object $X$ of $D^b(\md(A))$
determines an object of the stable module category $\stmod(A)$. For
example, if $X$ is the bounded complex
$$\dots\rightarrow X^{i-1}\rightarrow X^i\rightarrow
X^{i+1}\rightarrow\dots,$$
and if $X^i$ is projective for $i\neq0$,
then this object of $\stmod(A)$ is just $X^0$. A little more
generally, if all the terms of $X$ except for $X^n$ are projective,
then the corresponding object of $\stmod(A)$ is $\Omega^n(X^n)$. If
two objects of $D^b(\md(A))$ become isomorphic in $\stmod(A)$, we
shall say that they are \textbf{stably isomorphic}.

Since any equivalence of derived categories
$$D^b(\md(A))\approx D^b(\md(B))$$
restricts to an equivalence
$$K^b(\proj(A))\approx K^b(\proj(B)),$$
an equivalence of derived categories induces an equivalence of stable
module categories. This stable equivalence is `of Morita type'; i.e.,
it is induced (up to isomorphism) by an exact functor between the
module categories $\md(A)$ and $\md(B)$. 

Suppose we know that there is a stable equivalence of Morita type
$$F:\stmod(A)\approx\stmod(B),$$
and suppose we can produce an equivalence of derived categories
$$G:D^b(\md(B))\approx D^b(\md(C)),$$
where $C$ is a third
self-injective $k$-algebra: for example, we might construct a tilting
complex for $B$ and take $C$ to be its endomorphism algebra. Okuyama's
idea was to use a theorem~\cite[Theorem 2.1]{lin:Morita} of
Linckelmann, which states that if there is a stable equivalence of
Morita type
$$\stmod(A)\approx\stmod(C),$$
for self-injective algebras $A$ and
$C$, that takes the simple $A$-modules to the simple $C$-modules, then
$A$ and $C$ are Morita equivalent, and so their derived categories are
certainly equivalent. Thus if we can choose $G$ above so that the
induced equivalence between $\stmod(A)$ and $\stmod(C)$ has this property, then
$$D^b(\md(A))\approx D^b(\md(C))\approx D^b(\md(B)).$$

In a typical case of Brou\'{e}'s conjecture, where $A$ is a block
algebra of a group $G$ with abelian defect group $D$, and $B$ is the
Brauer correspondent block algebra of $N_G(D)$, the structure of $B$
is much easier to determine and understand than that of $A$. This is
the reason Okuyama's method is so useful for proving special cases of
Brou\'{e}'s conjecture: so long as we already know that there is a
stable equivalence of Morita type between $A$ and $B$, and so long as
we can determine the images of the simple modules under this
equivalence, the method requires no more information about $A$.

Okuyama has successfully used his method to verify several cases of
Brou\'{e}'s conjecture: for example, many cases involving defect group
$C_3\times C_3$.

\subsection{Combining Okuyama's method with Theorem~\ref{th:main}}
\label{sec:combine}

Since Theorem~\ref{th:main} not only provides an algebra $\Gamma$
whose derived category is equivalent to that of $\Lambda$, but also
identifies the objects that are sent to the simple $\Gamma$-modules,
it immediately combines with Okuyama's method to give the following
theorem.

\begin{thm}
  \label{th:equiv}
  Let $\Lambda$ and $\Gamma$ be finite-dimensional symmetric
  $k$-algebras, let 
  $$F:\md(\Gamma)\rightarrow\md(\Lambda)$$
  be an exact functor inducing a stable equivalence
  of Morita type, and let $\{S_0,\dots,S_r\}$ be a set of
  representatives for the isomorphism classes of simple
  $\Gamma$-modules.
  
  If there are objects $X_0,\dots,X_r$ of $D^b(\md(\Lambda)$ such
  that, for each $0\leq i\leq r$, $X_i$ is stably isomorphic to
  $F(S_i)$, and such that
  \begin{itemize}
  \item[(a)] $\Hom(X_i,X_j[m])=0$ for $m<0$,
  \item[(b)] $\Hom(X_i,X_j)=\Bigl\{
    \begin{array}{ll}
      0 & \textrm{if } i\neq j \\
      k & \textrm{if } i=j
    \end{array}$, and
  \item[(c)] $X_0,\dots,X_r$ generate $D^b(\md(\Lambda))$ as a
    triangulated category,
  \end{itemize}
  then $D(\Mod(\Lambda))$ and $D(\Mod(\Gamma))$ are equivalent as
  triangulated categories.
\end{thm}

\subsection{Objects with homology concentrated in one degree}
\label{sec:stalks}

As in Theorem~\ref{th:equiv}, let $\Lambda$ and $\Gamma$ be
finite-dimensional symmetric $k$-algebras,let
$$F:\md(\Gamma)\rightarrow\md(\Lambda)$$
be an exact functor inducing a
stable equivalence of Morita type, and let $\{S_0,\dots,S_r\}$ be a
set of representatives for the isomorphism classes of simple
$\Gamma$-modules. We shall assume that $F(S_i)$ is indecomposable for
every $i$: if this is not the case, then $F$ has a summand that
induces an isomorphic stable equivalence and which does send simple
modules to indecomposable modules, so this assumption does not involve
any real loss of generality.

In order to use Theorem~\ref{th:equiv} to prove that $\Lambda$ and
$\Gamma$ have equivalent derived categories, we need to find suitable
objects $X_0,\dots,X_r$ of $D^b(\md(\Lambda))$ that are stably
isomorphic to $F(S_0),\dots,F(S_r)$. We shall consider the case where
each of these objects has non-zero homology in only one degree. Since
$X_0,\dots,X_r$ must certainly be indecomposable if condition~(b) of
the theorem is to be satisfied, we must then have
\begin{equation}
\label{eq:syzygies}
X_i \cong \Omega^{n_i}F(S_i)[n_i]
\end{equation}
for some integers $n_0,\dots,n_r$.

Let us consider what is required in order for these objects to satisfy
conditions~(a) and (b) of Theorem~\ref{th:equiv}. It turns out that
many of the conditions are satisfied automatically, especially if
$n_0,\dots,n_r$ are small.

Recall that if $M$ and $N$ are modules for a finite-dimensional
self-injective $k$-algebra $A$, then for $m>0$ there are natural
isomorphisms
$$\Ext^m_A(M,N)\cong\Hom_{D(\Mod(A))}(M,N[m])\cong\sHom_A(\Omega^mM,N).$$

\begin{prop}
  \label{prop:stalks}
  Conditions~(a) and (b) of Theorem~\ref{th:equiv} are satisfied by
  the objects $X_i$ of (\ref{eq:syzygies}) if 
  \begin{itemize}
  \item[(i)] $\End_{\Lambda}(\Omega^{n_i}F(S_i))\cong k$ for $0\leq
    i\leq r$,
  \item[(ii)]
    $\Hom_{\Lambda}(\Omega^{n_i}F(S_i),\Omega^{n_j}F(S_j))=0$ whenever
    $i\neq j$ and $n_i\leq n_j$,
  \end{itemize}
  and if the four equivalent conditions
  \begin{itemize}
  \item[(iii)] $\sHom_{\Lambda}(\Omega^mF(S_i),F(S_j))=0$ whenever
    $-1>m>n_i-n_j$,
  \item[(iv)] $\sHom_{\Gamma}(\Omega^mS_i,S_j)=0$ whenever
    $-1>m>n_i-n_j$,
  \item[(v)] $\Ext^p_{\Lambda}(F(S_j),F(S_i))=0$ whenever
    $0<p<n_j-n_i-1$,
  \item[(vi)] $\Ext^p_{\Gamma}(S_j,S_i)=0$ whenever
    $0<p<n_j-n_i-1$,
  \end{itemize}
  are satisfied.
\end{prop}

Before we prove this, let us point out how few conditions this leaves
us to check when $n_0,\dots,n_r$ do not vary very much.

\begin{cor}
  \label{cor:stalks}
  Suppose
  $$\max_{0\leq i,j\leq r}|n_i-n_j| \leq 2.$$
  Then conditions~(a) and
  (b) of Theorem~\ref{th:equiv} are satisfied by the objects $X_i$ of
  (\ref{eq:syzygies}) if
  \begin{itemize}
  \item[(i)] $\End_{\Lambda}(\Omega^{n_i}F(S_i))\cong k$ for $0\leq
    i\leq r$, and
  \item[(ii)]
    $\Hom_{\Lambda}(\Omega^{n_i}F(S_i),\Omega^{n_j}F(S_j))=0$ whenever
    $i\neq j$ and $n_i\leq n_j$.
  \end{itemize}
\end{cor}

\begin{proof}[Proof of Proposition~\ref{prop:stalks}]
  The case $i=j$ of condition~(b) of Theorem~\ref{th:equiv} is just
  condition~(i) of the proposition.
  
  Since
  $$\Hom(X_i,X_j[m])=
  \Hom(\Omega^{n_i}F(S_i),\Omega^{n_j}F(S_j)[m+n_j-n_i]),$$
  we have
  $\Hom(X_i,X_j[m])=0$ when $m+n_j-n_i<0$, so we need only consider
  values of $m$ with $n_i-n_j\leq m\leq0$.
  
  The case $0\geq m=n_i-n_j$ requires
  $$0=\Hom_{\Lambda}(\Omega^{n_i}F(S_i),\Omega^{n_j}F(S_j))$$
  when
  $i\neq j$, which is condition~(ii) of the proposition.
  
  Since $F$ induces a stable equivalence,
  $$\sHom_{\Lambda}(\Omega^mF(S_i),F(S_j))\cong
  \sHom_{\Gamma}(\Omega^mS_i,S_j)$$
  for any $m$, and by Tate duality
  these spaces are dual to
  $$\sHom_{\Lambda}(\Omega^pF(S_j),F(S_i))\cong
  \sHom_{\Gamma}(\Omega^pS_j,S_i)$$
  for $p=-1-m$.
  
  What remains to be checked is that these spaces are zero for $0\geq
  m>n_i-n_j$, or $-1+n_j-n_i<p\leq-1$ (in which case $n_i\neq n_j$, so
  $i\neq j$). Since $S_i$ and $S_j$ are non-isomorphic simple modules,
  $$\sHom_{\Gamma}(S_i,S_j)=0=\sHom_{\Gamma}(S_j,S_i),$$
  so the cases
  $m=0$, $p=-1$ and $m=-1$, $p=0$ are automatically satisfied. The
  remaining cases are just conditions~(iii) to (vi) of the
  proposition.
\end{proof}

In fact, projective maps are the only obstruction to the automatic
satisfaction of condition~(i) and, if $n_j-n_i<2$, of condition~(ii).
For $\Lambda$-modules $M$ and $N$, write $\PHom_{\Lambda}(M,N)$ for
the space of projective maps from $M$ to $N$, so that
$$\sHom_{\Lambda}(M,N)=\Hom_{\Lambda}(M,N)/\PHom_{\Lambda}(M,N).$$

\begin{prop}
  \label{prop:projective}
  (a) Condition~(i) of Proposition~\ref{prop:stalks} (or
  Corollary~\ref{cor:stalks}) is satisfied for a given $i$ if and only
  if
  $$\PHom_{\Lambda}(\Omega^{n_i}F(S_i),\Omega^{n_i}F(S_i))=0.$$
  
  (b) Condition~(ii) is satisfied for a given $i$ and $j$ if and only
  if
  $$\PHom_{\Lambda}(\Omega^{n_i}F(S_i),\Omega^{n_j}F(S_j))=0$$
  and
  either
  \begin{itemize}
  \item[(I)] $n_j-n_i<2$, or
  \item[(II)] $n_j-n_i\geq2$ and
    $\Ext^{n_j-n_i-1}_{\Gamma}(S_j,S_i)=0.$
  \end{itemize}
\end{prop}

\begin{proof}
  Since $F$ induces a stable equivalence,
  $$\sHom_{\Lambda}(\Omega^{n_i}F(S_i),\Omega^{n_j}F(S_j))\cong
  \sHom_{\Gamma}(S_i,\Omega^{n_j-n_j}S_j),$$
  which is isomorphic to
  $k$ if $i=j$, is zero if $i\neq j$ and $n_i=n_j$, is dual to
  $\sHom_{\Gamma}(S_j,S_i)=0$ if $n_j-n_i=1$, and is dual to
  $\Ext^{n_j-n_i-1}_{\Gamma}(S_j,S_i)$ if $n_j-n_i>1$.
\end{proof}

\begin{rem}
Since there are no non-zero projective maps from a simple module to a
module with no projective summands, the vanishing of 
$\PHom_{\Lambda}(\Omega^{n_i}F(S_i),\Omega^{n_j}F(S_j))=0$
is automatic if $\Omega^{n_i}F(S_i)$ is a simple $\Gamma$-module.
\end{rem}

\section{Some examples}
\label{sec:examples}

We shall give some examples to show how Theorem~\ref{th:equiv} can be
used to verify specific cases of Brou\'{e}'s conjecture on
equivalences of derived categories for blocks with abelian defect
group.

\begin{conj}[Brou\'{e}]
  \label{conj:broue}
  Let $A$ be a block algebra of a finite group algebra $kG$, where $k$
  is an algebraically closed field of characteristic $p>0$. Suppose
  that a defect group $D$ of $A$ is abelian. Let $H=N_G(D)$, and let
  $B$ be the Brauer correspondent of $A$: a block algebra of
  $kH$. Then there is an equivalence $D(\Mod(A))\approx
  D(\Mod(B))$ of triangulated categories.
\end{conj}

We shall use the notation ($G$, $H$, $A$, $B$, $D$, $k$, $p$) of
Conjecture~\ref{conj:broue} in all the examples. In each case we shall
describe a set $\mathcal{X}=\{X_i:0\leq i\leq r\}$ of objects of
$D^b(\md(B))$ that are stably isomorphic to the images of the simple
$A$-modules under a stable equivalence. We shall then need to check
the conditions of Theorem~\ref{th:equiv}. The last of these conditions
is that $D^b(\md(B))$ is generated as a triangulated category by the
elements of $\mathcal{X}$. We shall use the notation
$\langle\mathcal{X}\rangle$ to refer to the triangulated category
generated by $\mathcal{X}$, and to prove that this is the whole of
$D^b(\md(B))$ we shall show that each of the simple $B$-modules is in
$\langle\mathcal{X}\rangle$. 

All of these examples are already known by other methods, but
Theorem~\ref{th:equiv} provides a simpler proof. In these cases we
shall give references to previous proofs.

In all the examples we give, it is well-known that there is a stable
equivalence which coincides with Green correspondence on objects. 

We start with a very simple example.

\subsection{Principal block of $G=A_5$, $p=2$}
\label{sec:a5}

The alternating group $G=A_5\cong SL(2,4)\cong PSL(2,5)$ has Sylow
2-subgroup $P\cong C_2\!\times\!C_2$, with normalizer $H=N_G(P)\cong
A_4\cong P\!\rtimes\!C_3$.

The principal block $A$ of $kG$ has three simple modules: the trivial
module and two 2-dimensional modules. 

$B=kH$ has three 1-dimensional simple modules, which we will denote by
$k$, $1$ and $2$.

The restrictions of the simple $A$-modules have the following
structures.
$$Y_0:=k,\;\;\;\;
Y_1:=
\begin{array}{c}
1 \\ 2 \\
\end{array},\;\;\;\;
Y_2:=
\begin{array}{c}
2 \\ 1 \\
\end{array}$$

We shall take $X_0:=k$ and $X_i:=\Omega Y_i[1]$ for
$i\in\{1,2\}$. The structure of $\Omega Y_i$ is as follows.
$$\Omega Y_1:=
\begin{array}{c}
k \\ 1 \\
\end{array},\;\;\;\;
\Omega Y_2:=
\begin{array}{c}
k \\ 2 \\
\end{array}$$

To verify conditions~(a) and (b) of Theorem~\ref{th:equiv}, it is
sufficient, by Corollary~\ref{cor:stalks}, to note that 
$$\End_B(k)\cong\End_B(\Omega Y_i)\cong k$$
and that 
$$\Hom_B(k,\Omega Y_i)=0$$
for $i\in\{1,2\}$.

Finally, to verify condition~(c) of Theorem~\ref{th:equiv}, note that
$k$ is certainly in $\langle\mathcal{X}\rangle$, and the other two
simples are in $\langle\mathcal{X}\rangle$ because they are the
kernels of surjective maps
$$\Omega Y_i\rightarrow k.$$

Hence $D^b(\md(A))\approx D^b(\md(B))$ by Theorem~\ref{th:equiv}.

The first published proof of this example is in~\cite{ric:splendid}.

\subsection{Principal block of $G=A_7$, $p=3$}
\label{sec:a7}

The alternating group $G=A_7$ has  Sylow
3-subgroup $P\cong C_3\!\times\!C_3$, with normalizer $H=N_G(P)\cong
P\!\rtimes\!C_4$.

The principal block $A$ of $kG$ has four simple modules: the trivial
module, two 10-dimensional modules and a 13-dimensional module.

$B=kH$ has four 1-dimensional simple modules, which we will denote by
$k$, $1$, $2$ and $3$. Fixing a generator $x$ for $C_4\leq H$ and a
primitive fourth root of unity $\zeta\in k$, $x$ acts on the simple
module $i$ as multiplication by $\zeta^i$.

The Green correspondents of the simple $A$-modules have the following
Loewy structures.
$$Y_0:=k,\;\;\;\;
Y_1:=1,\;\;\;\;
Y_2:=
\begin{array}{ccc}
&2& \\ 1&&3 \\ &2& \\
\end{array},\;\;\;\;
Y_3:=3$$

The Loewy structure of the projective cover of the simple module $2$ is:

$$P(2):=
\begin{array}{ccccc}
&&2&& \\ &1&&3& \\ k&&2&&k \\ &1&&3& \\ &&2&& \\
\end{array},$$
and so 
$$\Omega Y_2:=
\begin{array}{ccccc}
k&&&&k \\ &1&&3& \\ &&2&& \\
\end{array}.$$

We shall take $X_0:=k$, $X_1:=1$, $X_2:=\Omega Y_2[1]$, $X_3:=3$.

Since $\Omega Y_2$ has simple $2$ as its socle, and since this simple
occurs with multiplicity one as a composition factor of $\Omega Y_2$,
the conditions of Corollary~\ref{cor:stalks} are satisfied.

To check condition~(c) of Theorem~\ref{th:equiv}, note that the
simples $k$, $1$ and $3$ are certainly in $\langle\mathcal{X}\rangle$,
and hence $M:=\Omega Y_2/\soc(\Omega Y_2)$ is in
$\langle\mathcal{X}\rangle$, since each composition factor is
isomorphic to one of $k$, $1$, $3$. Therefore the simple $2$ is in
$\langle\mathcal{X}\rangle$, since it is the kernel of the natural
surjection $\Omega Y_2\rightarrow M$.

Okuyama~\cite{oku:derived} gave the first proof of this example.

\subsection{Principal block of $G=A_8$, $p=3$}
\label{sec:a8}

The alternating group $G=A_8$ has Sylow 3-subgroup $P\cong C_3\!\times\!C_3$,
with normalizer $H=N_G(P)\cong P\!\rtimes\!D_8$.

The principal block $A$ of $kG$ has five simple modules: the trivial
module and modules with dimensions 7, 13, 28 and 35.

$B=kH$ has four 1-dimensional modules, which we shall denote by $k$,
$1$, $2$ and $3$. We shall choose the names so that the kernel of the
action of $D_8$ on the simple $2$ is cyclic, whereas the kernels of
the actions on $1$ and $3$ are elementary abelian of rank two. There
is also a 2-dimensional simple module $S$.

The Green correspondents of the simple $A$-modules have the following
Loewy structures.

$$Y_0:=k,\;\;\;\;
Y_1:=2,\;\;\;\;
Y_2:=
\begin{array}{c}
1 \\ S \\ 1 \\
\end{array},\;\;\;\;
Y_3:=3,\;\;\;\;
Y_4:=S$$

We shall take $X_i:=Y_i$ for $i\in\{0,1,3,4\}$ and $X_2:=\Omega
Y_2[1]$.

The Loewy structure of the projective cover of the simple module $1$
is
$$\begin{array}{ccc}
&1& \\ &S& \\ k&1&2 \\ &S& \\ &1& \\
\end{array},$$
and so $\Omega Y_2$ has the structure
$$
\begin{array}{ccc}
k&&2 \\ &S& \\ &1& \\
\end{array}.$$

Since $\Omega Y_2$ has simple socle $1$, and since $1$ only occurs
with multiplicity one as a composition factor, $\End_B(\Omega
Y_2)\cong k$ and $\Hom_B(X_i,\Omega Y_2)=0$ for
$i\in\{0,1,3,4\}$. Hence the conditions of Corollary~\ref{cor:stalks}
are satisfied. 

Condition~(c) of Theorem~\ref{th:equiv} also holds, because clearly
all the simples other than $1$ are in $\langle\mathcal{X}\rangle$, and
since $\Omega Y_2$, which contains $1$ as a composition factor with
multiplicity one, is also in $\langle\mathcal{X}\rangle$, it follows
that $1$ is in $\langle\mathcal{X}\rangle$.

Okuyama~\cite{oku:derived} gave the first proof of this example.

\section{General coefficient fields}
\label{sec:non-closed}

The only use we made of the condition that $k$ is algebraically closed
was to assume that the endomorphism ring of a simple module for a
finite-dimensional $k$-algebra was just $k$. Our main theorem,
Theorem~\ref{th:main}, is true a little more generally for a field $k$
that is not algebraically closed, as we can allow $\End(X_i)$ to be
any finite-dimensional division algebra over $k$. The only difference
in the construction of the objects $T_i$ is that rather than using a
$k$-basis $B_i^{(n-1)}(j,t)$ of $\Hom(X_j[t],X_i^{(n-1)})$ in order to
form the object $Z_i^{(n-1)}(j,t)$, we should use a basis as a left
$\End(X_j)$-module.

Theorem~\ref{th:equiv} generalizes similarly.


\begin{thebibliography}{10}

\bibitem{BN:hocolim}
Marcel B{\"o}kstedt and Amnon Neeman.
\newblock Homotopy limits in triangulated categories.
\newblock {\em Compositio Math.}, 86(2):209--234, 1993.

\bibitem{bro:perfect}
Michel Brou{\'e}.
\newblock Isom\'etries parfaites, types de blocs, cat\'egories d\'eriv\'ees.
\newblock {\em Ast\'erisque}, (181-182):61--92, 1990.

\bibitem{chu:sl2}
Joseph Chuang.
\newblock Derived equivalence in ${S}{L}\sb 2(p\sp 2)$.
\newblock {\em Trans. Amer. Math. Soc.}, 353(7):2897--2913 (electronic), 2001.

\bibitem{hol:thesis}
Miles Holloway.
\newblock {\em Derived equivalences for group algebras}.
\newblock PhD thesis, Bristol, 2001.

\bibitem{lin:Morita}
Markus Linckelmann.
\newblock Stable equivalences of {M}orita type for self-injective algebras and
  $p$-groups.
\newblock {\em Math. Z.}, 223(1):87--100, 1996.

\bibitem{oku:derived}
Tetsuro Okuyama.
\newblock Some examples of derived equivalent blocks of finite groups.
\newblock Preprint, Hokkaido, 1998.

\bibitem{oku:sl2}
Tetsuro Okuyama.
\newblock Derived equivalence in ${S}{L}(2,q)$.
\newblock Preprint, Hokkaido, 2000.

\bibitem{ric:Morita}
Jeremy Rickard.
\newblock Morita theory for derived categories.
\newblock {\em J. London Math. Soc. (2)}, 39(3):436--456, 1989.

\bibitem{ric:equivalences}
Jeremy Rickard.
\newblock Derived equivalences as derived functors.
\newblock {\em J. London Math. Soc. (2)}, 43(1):37--48, 1991.

\bibitem{ric:splendid}
Jeremy Rickard.
\newblock Splendid equivalences: derived categories and permutation modules.
\newblock {\em Proc. London Math. Soc. (3)}, 72(2):331--358, 1996.

\bibitem{rou:gluing}
Rapha\"{e}l Rouquier.
\newblock Gluing $p$-permutation modules.
\newblock Preprint, Paris, 1988.

\bibitem{rou:block}
Rapha{\"e}l Rouquier.
\newblock Block theory via stable and {R}ickard equivalences.
\newblock In {\em Modular representation theory of finite groups}, pages
  101--146. de Gruyter, Berlin, 2001.

\end{thebibliography}
\bibliographystyle{plain}

\end{document}